\documentclass[12pt]{amsart}
\usepackage{amsmath,amsthm,latexsym,amscd,amsbsy,amssymb,url}
\usepackage[all]{xypic}
 \relax

\chardef\bslash=`\\ 

\makeatletter
\def\verbatim{\interlinepenalty\@M \@verbatim
  \leftskip\@totalleftmargin\advance\leftskip2pc
  \frenchspacing\@vobeyspaces \@xverbatim}
\makeatother
\makeatletter
  \def\dgt@k{\dg@DX=-3 \dg@DY=2 \dg@SIZE=3}
\makeatother

\makeatletter
  \def\dgt@kk{\dg@DX=3 \dg@DY=-1 \dg@SIZE=3}%
\makeatother

\theoremstyle{plain}
\newtheorem{thm}{Theorem}[section]
\newtheorem{cor}[thm]{Corollary}
\newtheorem{lem}[thm]{Lemma}

\theoremstyle{definition}

\numberwithin{equation}{section}

\begin{document}


\title
{Retracts of sigma-products of Hilbert cubes}
\author{A. Chigogidze}
\address{Department of Mathematics and Statistics,
University of North Carolina at Greensboro,
383 Bryan Bldg, Greensboro, NC, 27402, USA}
\email{chigogidze@www.uncg.edu}


\begin{abstract}{We consider the sigma-product of the $\omega_{1}$-power of the Hilbert cube. This space is characterized among its retracts as the only one without $G_{\delta}$-points.}
\end{abstract}

\maketitle
\markboth{A.~Chigogidze}
{Topology of $\Sigma$-products of $\omega_{1}$-many Hilbert cubes}

\section{Introduction}
Sigma-products and their subspaces have been extensively studied by topologists and functional analysts for several decades. We refer the reader to \cite{ka} where a comprehensive survey of the related results from both topology and functional analysis are discussed in detail.

Recall that the sigma-product $\Sigma (X, \ast )$ of an uncountable collection $\{ X_{t} \colon t \in T\}$ of spaces with base points $\ast_{t} \in X_{t}$, $t \in T$, is the subspace of the product $X= \prod\{ X_{t} \colon t \in T \}$ defined as follows

\[ \Sigma (X, \ast ) = \{ \{ x_{t} \colon t \in T\} \in X \colon |\{ t \in T \colon x_{t}\neq \ast_{t}|\leq \omega\}\} .\]

We are interested in the case when each $X_{t}$ is a copy of the Hilbert cube $I^{\omega}$, $|T| = \omega_{1}$ and $\ast_{t}$ is the point (in $I^{\omega}$) all coordinates of which equal to $0$. The corresponding sigma-product is denoted by $\Sigma$.

Our main result (Theorem \ref{T:characterization}) states that if a retract of $\Sigma$ does not contain $G_{\delta}$-points, then it is homeomorphic to $\Sigma$.


\section{Auxiliary lemmas}
Terminology, notation and results related to inverse spectra and absolute retracts used here can be found in \cite{chibook}. One of the main concepts we need below is that of $\omega$-spectra ${\mathcal S} = \{ X_{\alpha}, p_{\alpha}^{\beta}, A\}$. These are $\omega$-continuous inverse spectra consisting of metrizable compact spaces $X_{\alpha}$, surjective projections $p_{\alpha}^{\beta} \colon X_{\beta} \to X_{\alpha}$, $\alpha \leq \beta$, and an $\omega$-complete indexing set $A$. This essentially means that $A$ contains supremums of countable chains and that for any such chain $\{ \alpha_{n} \colon n\in \omega\}$ the space $X_{\alpha}$, where $\alpha = \sup\{ \alpha_{n} \colon n\in \omega\}$, is naturally homeomorphic to the limit of the inverse sequence ${\mathcal S}_{\alpha} = \{ X_{\alpha_{n}}, p_{\alpha_{n}}^{\alpha_{n+1}},\omega\}$.

Recall also that a compact space is an absolute retract if and only if it is a retract of a Tychonov cube and that a map $p \colon X \to Y$ of compact spaces is soft if for any compactum $B$, 
      it's closed subset $A$ and maps $g$ and $h$ such that the following diagram (of undotted arrows)
      commutes
        \[
        \xymatrix{
            X \ar^{p}[r] & Y \\
            A \ar@^{(->}^{i}[r] \ar^{g}[u] & B \ar_{h}[u]  \ar@^{.>}_{
k}[ul]\\
        }
      \]

\noindent there exists a map $k \colon B \to X$ (the dotted arrow) such that $k |A = g$ and $fk = h$.

The prime example of a soft map is the projection $X \times I^{\omega} \to X$.

Finally recall that for a given map $p \colon X \to Y$ a closed subset $F \subseteq X$ is a fibered $Z$-set in $X$ (with respect to $p$) if for any open cover ${\mathcal U} \in \operatorname{cov}(X)$ there exists a map $f_{\mathcal U} \colon X\to X$ such that $pf_{\mathcal U} = p$ (i.e. $f_{\mathcal U}$ acts fiberwise), $f_{\mathcal U}(X) \cap F = \emptyset$ and $f_{\mathcal U}$ is ${\mathcal U}$-close to $\operatorname{id}_{X}$. 

\begin{lem}\label{L:disjointsection}
Let a non-metrizable compact space $X$ be represented as the limit of an $\omega$-spectrum ${\mathcal S} = \{ X_{\alpha}, p_{\alpha}^{\beta}, A\}$ with soft projections $p_{\alpha}^{\beta}$. Suppose that $F$ is a closed subset of $X$ containing no closed $G_{\delta}$-subsets of $X$. Then for each $\alpha \in A$ there exists $\beta \in A$, with $\beta > \alpha$, such that there is a map $i_{\alpha}^{\beta} \colon X_{\alpha} \to X_{\beta}$ satisfying the following two properties:
\begin{itemize}
\item[(1)]
$p_{\alpha}^{\beta}i_{\alpha}^{\beta} = \operatorname{id}_{X_{\alpha}}$,
\item[(2)]
$i_{\alpha}^{\beta}(X_{\alpha}) \cap p_{\beta}(F) = \emptyset$.
\end{itemize}
\end{lem}
\begin{proof}
Let $\alpha_{0} = \alpha$ and $x_{0} \in X_{\alpha_{0}}$. Since $p_{\alpha_{0}}^{-1}(x_{0})$ is closed and $G_{\delta}$ in $X$, it follows from the assumption that $p_{\alpha_{0}}^{-1}(x_{0}) \setminus F \neq \emptyset$. Take an index $\alpha_{1}\in A$ such that $\alpha_{1}>\alpha_{0}$ and $(p_{\alpha_{0}}^{\alpha_{1}})^{-1}(x_{0}) \setminus p_{\alpha_{1}}(F) \neq \emptyset$. Let $x_{1} \in (p_{\alpha_{0}}^{\alpha_{1}})^{-1}(x_{0}) \setminus p_{\alpha_{1}}(F)$. The softness of the projection $p_{\alpha_{0}}^{\alpha_{1}} \colon X_{\alpha_{1}} \to X_{\alpha_{0}}$ guarantees the existence of a map $i_{0}^{1} \colon X_{\alpha_{0}} \to X_{\alpha_{1}}$ such that
$p_{\alpha_{0}}^{\alpha_{1}}i_{0}^{1} = \operatorname{id}_{X_{\alpha_{0}}}$ and $i_{0}^{1}(x_{0}) = x_{1}$. Let 
\[ V_{1} = \{ x \in X_{\alpha_{0}}\colon i_{0}^{1}(x) \notin p_{\alpha_{1}}(F)\}.\]

Note that $x_{0} \in V_{1}$ and consequently $V_{1}$ is a non-empty open subset of $X_{\alpha_{0}}$.

Let $\gamma < \omega_{1}$. Suppose that for each $\lambda$, $1 \leq \lambda < \gamma$, we have already constucted an index $\alpha_{\lambda} \in A$, an open subset $V_{\lambda} \subseteq X_{\alpha_{0}}$ and a section $i_{0}^{\lambda} \colon X_{\alpha_{0}} \to X_{\alpha_{\lambda}}$ of the projection $p_{\alpha_{0}}^{\alpha_{\lambda}} \colon X_{\alpha_{\lambda}} \to X_{\alpha_{0}}$, satisfying the following conditions:
\begin{itemize}
\item[(i)]
$\alpha_{\lambda} < \alpha_{\mu}$, wherenever $\lambda < \mu < \gamma$,
\item[(ii)]
$\alpha_{\mu} = \sup\{ \alpha_{\lambda} \colon \lambda < \mu\}$, whenever $\mu < \gamma$ is a limit ordinal,
\item[(iii)]
$V_{\lambda} \subsetneq V_{\mu}$, whenever $\lambda < \mu < \gamma$,
\item[(iv)]
$V_{\mu} = \cup\{ V_{\lambda} \colon \lambda < \mu\}$, whenever $\mu < \gamma$ is a limit ordinal,
\item[(v)]
$i_{0}^{\mu} = \triangle\{ i_{0}^{\lambda} \colon \lambda < \mu\}$, whenever $\mu < \gamma$ is a limit ordinal,
\item[(vi)]
$i_{0}^{\lambda} = p_{\alpha_{\lambda}}^{\alpha_{\mu}}i_{0}^{\mu}$, whenever $\lambda < \mu < \gamma$,
\item[(vii)]
$V_{\lambda} = \{ x \in X_{\alpha_{0}}\colon i_{0}^{\lambda}(x) \notin p_{\alpha_{\lambda}}(F)\}$.
\end{itemize}
We shall construct the index $\alpha_{\gamma}$, the open subset $V_{\gamma} \subseteq X_{\alpha_{0}}$ and the section $i_{0}^{\gamma} \colon X_{\alpha_{0}} \to X_{\alpha_{\gamma}}$ of the projection $p^{\alpha_{\gamma}}_{\alpha_{0}} \colon X_{\alpha_{\gamma}} \to X_{\alpha_{0}}$.

Suppose that $\gamma$ is a limit ordinal. By (i), $\{ \alpha_{\mu} \colon \mu < \gamma\}$ is a countable chain in $A$ and we let (recall that the indexing set $A$ of ${\mathcal S}$ is a $\omega$-compelete set and therefore contains supremums of countable chains of its elements)
\[ \alpha_{\gamma} = \sup\{ \alpha_{\mu} \colon \mu < \gamma\} \in A .\]

By the $\omega$-continuity of the spectrum ${\mathcal S}$, the compactum $X_{\alpha_{\gamma}}$ is naturally homeomorphic to the limit of the inverse sequence $\{ X_{\alpha_{\mu}}, p_{\alpha_{\lambda}}^{\alpha_{\mu}}, \lambda, \mu < \gamma\}$. Consequently, by (vi), the diagonal product
\[ i_{0}^{\gamma} = \triangle\{ i_{0}^{\mu} \colon \mu < \gamma\} \colon X_{\alpha_{0}}\to X_{\alpha_{\gamma}}\]

\noindent is well-defined and satifies corresponding conditions (v) and (vi). Let 
\[ V_{\gamma} = \{ x \in X_{\alpha_{0}}\colon i_{0}^{\gamma}(x) \notin p_{\alpha_{\gamma}}(F)\} .\]

\noindent Note that $V_{\gamma} = \cup\{ V_{\alpha_{\mu}} \colon \mu < \gamma\}$. Then, corresponding conditions (vii), (iii) and (iv) are satisfied.

Next consider the case $\gamma = \mu +1$. In case $V_{\mu} = X_{\alpha_{0}}$, the desired $\beta$ is $\alpha_{\mu}$. Suppose that $V_{\mu} \neq X_{\alpha_{0}}$ and let
\[ x_{\mu} = i_{0}^{\mu}(z) \in i_{0}^{\mu}(X_{\alpha_{0}}) \subseteq X_{\alpha_{\mu}},\]

\noindent where $z \in X_{\alpha_{0}}\setminus V_{\mu}$. Since $p_{\alpha_{\mu}}^{-1}(x_{\mu})$ is closed and $G_{\delta}$ in $X$, we have $p_{\alpha_{\mu}}^{-1}(x_{\mu}) \setminus F \neq \emptyset$ (note that $X_{\alpha_{\mu}}$ is a metrizable compactum). Choose an index $\alpha_{\gamma} \in A$ so that $\alpha_{\gamma} > \alpha_{\mu}$ and $(p_{\alpha_{\mu}}^{\alpha_{\gamma}})^{-1}(x_{\mu}) \setminus p_{\alpha_{\gamma}}(F) \neq \emptyset$.

Softness of the projection $p_{\alpha_{\mu}}^{\alpha_{\gamma}} \colon X_{\alpha_{\gamma}} \to X_{\alpha_{\mu}}$ guarantees the existence of a map $i_{\mu}^{\gamma} \colon X_{\alpha_{\mu}} \to X_{\alpha_{\gamma}}$ such that $p_{\alpha_{\mu}}^{\alpha_{\gamma}}i_{\mu}^{\gamma} = \operatorname{id}_{X_{\alpha_{\mu}}}$ and $i_{\mu}(x_{\mu}) = z^{\prime}$, where $z^{\prime} \in (p_{\alpha_{\mu}}^{\alpha_{\gamma}})^{-1}(x_{\mu}) \setminus p_{\alpha_{\gamma}}(F)$. Let $i_{0}^{\gamma} = i_{\mu}^{\gamma}i_{0}^{\mu}$ and $V_{\gamma} = \{ x \in X_{\alpha_{0}}\colon i_{0}^{\gamma}(x) \notin p_{\alpha_{\gamma}}(F)\}$. Note that $V_{\mu} \subseteq V_{\gamma}$ and $z^{\prime} \in V_{\gamma}\setminus V_{\mu}$. This completes construction of the needed objects in the case $\gamma = \mu +1$.

Thus the construction can be carried out for each $\lambda < \omega_{1}$ and we obtain a strictrly increasing collection $\{ V_{\lambda} \colon \lambda < \omega_{1}\}$ of open subsets of the metrizable compactum $X_{\alpha_{0}}$. Clearly, this collection must stabilize, which means that there in an index $\lambda_{0} < \omega_{1}$ such that $V_{\lambda} = V_{\lambda_{0}}$ for any $\lambda \geq \lambda_{0}$.
By construction, this is only possible if $V_{\lambda_{0}} = X_{\alpha_{0}}$. Let $\beta = \alpha_{\lambda_{0}}$ and $i_{\alpha}^{\beta} = i_{\alpha_{0}}^{\lambda_{0}}$. Clearly $i_{\alpha}^{\beta}(X_{\alpha}) \cap p_{\beta}(F) = \emptyset$.
\end{proof}

We also need the following statement.

\begin{lem}\label{L:fiberedzset}
Let a non-metrizable compact space $X$ be represented as the limit of an $\omega$-spectrum ${\mathcal S} = \{ X_{\alpha}, p_{\alpha}^{\beta}, A\}$ with soft projections $p_{\alpha}^{\beta}$. Suppose that $F$ is a closed subset of $X$ containing no closed $G_{\delta}$-subsets of $X$. Then for each $\alpha \in A$ there exists an index $\beta \in A$, with $\beta > \alpha$, such that $p_{\beta}(F)$ is a fibered $Z$-set in $X_{\beta}$ with respect to the projection $p_{\alpha}^{\beta} \colon X_{\beta} \to X_{\alpha}$.
\end{lem}
\begin{proof}
Let $\alpha_{0} = \alpha$. Choose $\alpha_{k+1} > \alpha_{k}$ so that the projection $p_{\alpha_{k}}^{\alpha_{k+1}} \colon X_{\alpha_{k+1}} \to X_{\alpha_{k}}$ has a section $i_{k}^{k+1} \colon X_{\alpha_{k}} \to X_{\alpha_{k+1}}$ such that $i_{k}^{k+1}(X_{\alpha_{k}}) \cap p_{\alpha_{k+1}}(F) = \emptyset$. Let $\beta = \sup\{ \alpha_{k} \colon k \in \omega\}$. Let us show that $p_{\beta}(F)$ is a fibered $Z$-set in $X_{\beta}$ with respect to the projection $p^{\beta}_{\alpha}$. Let 
${\mathcal U} = \{ U_{i} \colon i \in I\}$ be an open cover of $X_{\beta}$. Without loss of generality we may assume that $U_{i} = (p_{\alpha_{k}}^{\beta})^{-1}(U_{i}^{k})$, $i \in I$, where $k \in \omega$ and $U_{i}^{k}$ is open in $X_{\alpha_{k}}$. Let $j \colon X_{\alpha_{k+1}}\to X_{\beta}$ be any section of the projection $p_{\alpha_{k+1}}^{\beta} \colon X_{\beta} \to X_{\alpha_{k+1}}$. Consider the map $f_{\mathcal U} = ji_{k}^{k+1}p_{\alpha_{k}}^{\beta} \colon X_{\beta} \to X_{\beta}$. Since
\[ p_{\alpha_{k}}^{\beta}f_{\mathcal U} = p_{\alpha_{k}}^{\beta}ji_{k}^{k+1}p_{\alpha_{k}}^{\beta} = p_{\alpha_{k}}^{\alpha_{k+1}}(p_{\alpha_{k+1}}^{\beta}j)i_{k}^{k+1}p_{\alpha_{k}}^{\beta} = (p_{\alpha_{k}}^{\alpha_{k+1}}i_{k}^{k+1})p_{\alpha_{k}}^{\beta} = p_{\alpha_{k}}^{\beta},\]
\noindent it follows that $f_{\mathcal U}$ is ${\mathcal U}$-close to $\operatorname{id}_{X_{\beta}}$. Also $p_{\alpha}^{\beta}f_{\mathcal U} = p_{\alpha}^{\alpha_{k}}p_{\alpha_{k}}^{\beta}f_{\mathcal U}=p_{\alpha}^{\alpha_{k}}p_{\alpha_{k}}^{\beta} = p_{\alpha}^{\beta}$ (i.e. $f_{\mathcal U}$ acts fiberwise with respect to $p_{\alpha}^{\beta}$). It only remains to note $f_{\mathcal U}(X_{\beta}) \cap p_{\beta}(F) = \emptyset$.
\end{proof}

\begin{lem}\label{L:pseudocompact}
Let $X$ be a pseudocompact space without $G_{\delta}$-points. If $\beta X$ - its Stone-\v{C}ech compactification - is an absolute retract of weight $\omega_{1}$, then
$\beta X$ is homeomorphic to $I^{\omega_{1}}$.
\end{lem}
\begin{proof}
By assumption, if $x \in X$, then $x$ is not a $G_{\delta}$-point in $\beta X$. Since $X$ is pseudocompact, no point in $\beta X \setminus X$ is a $G_{\delta}$-subset in $\beta X$ (see \cite[Exercise 6I.1]{gj}). Thus $\beta X$ has no $G_{\delta}$-points. By \v{S}\v{c}epin's theorem (see \cite[Theorem 7.2.9]{chibook}), $\beta X \approx I^{\omega_{1}}$.
\end{proof}

\section{Main result}
In this section we prove our main result.
\begin{thm}\label{T:characterization}
Let $X$ be a retract of $\Sigma$. Then the following conditions are equivalent:
\begin{itemize}
\item[(i)]
$X$ is homeomorphic to $\Sigma$,
\item[(ii)]
$X$ has no $G_{\delta}$-points.
\end{itemize}
\end{thm}
\begin{proof}
The implication (i) $\Rightarrow$ (ii) is trivial.

(ii) $\Rightarrow$ (i). Let $|A| = \omega_{1}$. First let us introduce some notation. If $C \subset B \subseteq A$, then $\pi_{B} \colon (I^{\omega})^{A} \to (I^{\omega})^{B}$ and $\pi_{C}^{B} \colon (I^{\omega})^{B} \to (I^{\omega})^{C}$ denote the corresponding projections. Similarly by $\lambda_{B} \colon (I^{\omega})^{B} \to (I^{\omega})^{A}$ and $\lambda_{C}^{B} \colon (I^{\omega})^{C} \to (I^{\omega})^{B}$ we denote the sections of $\pi_{B}$ and $\pi_{C}^{B}$ defined as follows:
\[ \lambda_{B}(\{ x_{t} \colon t \in B\}) = \left(\{ x_{t} \colon t \in B\} ,\{ 0_{t} \colon t \in A\setminus B\}\right) \]

\noindent and

\[ \lambda_{C}^{B}(\{ x_{t} \colon t \in C\}) = \left(\{ x_{t} \colon t \in C\} ,\{ 0_{t} \colon t \in B\setminus C\}\right) .\]

\noindent Here $0_{t}$ denotes the point in the $t$-th copy of the Hilbert cube, all coordinates of which ere equal to $0$. Note that $\Sigma = \bigcup\{ \lambda_{B}((I^{\omega})^{B}) \colon B \in \exp_{\omega}A\}$.

Let $X \subseteq \Sigma$ and $r \colon \Sigma \to X$ be a retraction. Recall that $\Sigma$ is normal, pseudocompact and $\beta \Sigma = (I^{\omega})^{A}$ (see \cite[Problems 2.7.14, 3.10.E and 3.12.23(c)]{eng}). Consequently, $\beta X = \operatorname{cl}_{(I^{\omega})^{A}}X$ and $r$ has the extension $\tilde{r} \colon (I^{\omega})^{A} \to \operatorname{cl}_{(I^{\omega})^{A}}X =Y$. Note that $\tilde{r}$ is also a retraction and consequently $Y$ is a compact absolute retract. Note that $X$, as a retract of $\Sigma$, is pseudocompact. Therefore, by Lemma \ref{L:pseudocompact}, $Y \approx I^{\omega_{1}}$.  

For each $C, B \subseteq A$, with $C \subseteq B$, let $Y_{B} = \pi_{B}(Y)$, $p_{B} = \pi_{B}|Y$ and $p_{C}^{B} = \pi_{C}^{B}|Y_{B}$. Clearly, $Y$ is the limit space of the $\omega$-spectrum ${\mathcal S}_{Y} = \{ Y_{B}, p_{C}^{B}, \exp_{\omega}A\}$. Since $Y \approx I^{\omega_{1}}$, \v{S}\v{c}epin's spectral theorem (see \cite[Theorem 1.3.4]{chibook}) for $\omega$-spectra insures that there exists an $\omega$-closed and cofinal subset ${\mathcal A} \subseteq \exp_{\omega}A$ such that 
\begin{enumerate}
\item
$Y = \lim{\mathcal S}_{\mathcal A}$, where ${\mathcal S}_{\mathcal A} = \{ Y_{B}, p_{C}^{B}, {\mathcal A}\}$,
\item
$Y_{B} \approx I^{\omega}$ whenever $B \in {\mathcal A}$,
\item
$p_{C}^{B} \colon Y_{B} \to Y_{C}$ is a trivial fibration with fiber $I^{\omega}$, whenever $C \subseteq B$, $C, B \in {\mathcal A}$.
\end{enumerate}

Applying the same spectral theorem to the map $\tilde{r} \colon (I^{\omega})^{A} \to Y$ and to the $\omega$-spectra ${\mathcal S} = \{ (I^{\omega})^{B}, \pi_{C}^{B}, {\mathcal A}\}$ (whose limit is $(I^{\omega})^{A}$) and ${\mathcal S}_{\mathcal A}$ we can find an $\omega$-closed and cofinal subset ${\mathcal B} \subseteq \exp_{\omega}A$ such that ${\mathcal B} \subseteq {\mathcal A}$ and for each $B \in {\mathcal B}$ there exists a retraction $r_{B} \colon (I^{\omega})^{B} \to Y_{B}$ such that $p_{B}\tilde{r} = r_{B}\pi_{B}$.

Let $B \in {\mathcal B}$ and consider the composition $i_{B} = \tilde{r}\lambda_{B}|Y_{B} \colon Y_{B} \to Y$. Note that $p_{B}i_{B} = p_{B}\tilde{r}\lambda_{B}|Y_{B} = r_{B}\pi_{B}\lambda_{B}|Y_{B} = r_{B}|Y_{B} = \operatorname{id}_{Y_{B}}$. In other words, $i_{B}$ is a section of the projection $p_{B}$. If $C \subseteq B$, $C, B \in {\mathcal B}$, we let $i_{C}^{B} = p_{B}i_{C}$. Note that $i_{C}^{B}$ is a section of the projection $p_{C}^{B}$.

Next we show that $X = \bigcup\{ i_{B}(Y_{B}) \colon B \in {\mathcal B}\}$. Indeed, let $x \in X$. Since $\Sigma = \bigcup\{ \lambda_{B}((I^{\omega})^{B})\colon B \in {\mathcal B}\}$, there is $B \in {\mathcal B}$ such that $x = \lambda_{B}(y)$ for some $y \in (I^{\omega})^{B}$. But $y = \pi_{B}(\lambda_{B}(y)) = \pi_{B}(x) \subseteq \pi_{B}(X) \subseteq \pi_{B}(Y)= Y_{B}$. Consequently, $i_{B}(y) = 
\tilde{r}(\lambda_{B}(y)) =\tilde{r}(x) = r(x) = x$.

Next we will construct a cofinal collection of countable subsets $\{A_{\alpha} \colon \alpha < \omega_{1}\} \subseteq {\mathcal B}$ of $A$ and homeomorphisms $h_{\alpha} \colon Y_{A_{\alpha}} \to (I^{\omega})^{A_{\alpha}}$ satisfying the following conditions:
\begin{itemize}
\item[(i)]
$A_{\alpha} \subseteq A_{\beta}$, whenever $\alpha < \beta < \omega_{1}$;
\item[(ii)]
$A_{\beta} = \bigcup\{ A_{\alpha} \colon \alpha < \beta\}$, whenever $\beta <\omega_{1}$ is a limit ordinal;
\item[(iii)]
For each $\alpha < \omega_{1}$, $\pi_{A_{\alpha}}^{A_{\alpha+1}}h_{\alpha +1} = h_{\alpha}p_{A_{\alpha}}^{A_{\alpha+1}}$, i.e. the following diagram is commutative
 \[
        \xymatrix{
            Y_{A_{\alpha +1}} \ar^{h_{\alpha +1}}[r] \ar_{p_{A_{\alpha}}^{A_{\alpha +1}}}[d]  &(I^{\omega})^{A_{\alpha +1}} \ar^{\pi_{A_{\alpha}}^{A_{\alpha +1}}}[d]\\
            Y_{A_{\alpha}} \ar^{h_{\alpha}}[r] &  (I^{\omega})^{A_{\alpha}}   \\
        }
      \]
 
\item[(iv)]
$h_{\beta} = \lim\{ h_{\alpha} \colon \alpha < \beta\}$, whenever $\beta <\omega_{1}$ is a limit ordinal;
\item[(v)]
$p_{A_{\alpha}}^{A_{\alpha +1}} \colon Y_{A_{\alpha +1}}\to Y_{A_{\alpha}}$ is a trivial fibration with fiber $I^{\omega}$;
\item[(vi)]
$i_{A_{\alpha}}^{A_{\alpha+1}}(Y_{A_{\alpha}})$ is a fibered $Z$-set in $Y_{A_{\alpha +1}}$ with respect to the projection $p_{A_{\alpha}}^{A_{\alpha +1}}$;
\item[(vii)]
For each $\alpha < \omega_{1}$, $h_{\alpha +1}i_{A_{\alpha}}^{A_{\alpha+1}} = \lambda_{A_{\alpha}}^{A_{\alpha +1}}h_{\alpha}$, i.e. the following diagram commutes:

\[
        \xymatrix{
            Y_{A_{\alpha +1}} \ar^{h_{\alpha +1}}[r]   &(I^{\omega})^{A_{\alpha +1}} \\
            Y_{A_{\alpha}} \ar^{h_{\alpha}}[r] \ar^{i_{A_{\alpha}}^{A_{\alpha +1}}}[u]&  (I^{\omega})^{A_{\alpha}}\;.   \ar_{\lambda_{A_{\alpha}}^{A_{\alpha +1}}}[u] \\
        }
      \]
\end{itemize}

Let $A_{0}$ be any element of ${\mathcal B}$ and take any homeomorphism $h_{0} \colon Y_{A_{0}}\to (I^{\omega})^{A_{0}}$. 

Let $\beta < \omega_{1}$. Suppose that for each $\alpha < \beta$ we have already constructed a countable set $A_{\alpha} \in {\mathcal B}$ and a homeomorphism $h_{\alpha} \colon Y_{A_{\alpha}} \to (I^{\omega})^{A_{\alpha}}$ satisfying the above conditions for appropriate indices. We proceed by constructing these objects for the ordinal $\beta$.

If $\beta = \sup\{ \alpha \colon \alpha < \beta\}$, then set $A_{\beta} = \cup\{ A_{\alpha} \colon \alpha < \beta\}$ and $h_{\beta} = \lim\{ h_{\alpha} \colon \alpha < \beta\}$. Then, all required conditions are clearly satisfied.

Now consider the case $\beta = \alpha +1$. Since $i_{A_{\alpha}}(Y_{A_{\alpha}})$ is a metrizable compactum in $Y$, it cannot contain closed $G_{\delta}$-subsets of $Y$ (which contain copies of the Tychonov cube $I^{\omega_{1}}$). Consequently, we can find, based on Lemma \ref{L:fiberedzset}, an element $A_{\alpha +1} \in {\mathcal B}$, such that $A_{\alpha} \subseteq A_{\alpha +1}$ and $i_{A_{\alpha}}^{A_{\alpha +1}}(Y_{A_{\alpha}}) = p_{A_{\alpha +1}}(i_{A_{\alpha}}(Y_{A_{\alpha}}))$ is a fibered $Z$-set with respect to the projection $p_{A_{\alpha}}^{A_{\alpha +1}}$. Since both projections $p_{A_{\alpha}}^{A_{\alpha +1}}$ and $\pi_{A_{\alpha}}^{A_{\alpha +1}}$ are trivial fibrations with fiber $I^{\omega}$, there exists a homeomorphism $f \colon Y_{A_{\alpha +1}} \to (I^{\omega})^{A_{\alpha +1}}$ such that $\pi_{A_{\alpha}}^{A_{\alpha +1}}f = h_{A_{\alpha}}p_{A_{\alpha}}^{A_{\alpha +1}}$. Then the set $f(i_{A_{\alpha}}^{A_{\alpha +1}}(Y_{A_{\alpha}}))$ is a fibered $Z$-set in $(I^{\omega})^{A_{\alpha +1}}$ with respect to the projection $\pi_{A_{\alpha}}^{A_{\alpha +1}}$. Consider now another fibered $Z$-set in $(I^{\omega})^{A_{\alpha +1}}$ (also with respect to the projection $\pi_{A_{\alpha}}^{A_{\alpha +1}}$) -- namely, $\lambda_{A_{\alpha}}^{A_{\alpha +1}}((I^{\omega})^{A_{\alpha}})$. There is a homeomorphism 

\[ g \colon f(i_{A_{\alpha}}^{A_{\alpha +1}}(Y_{A_{\alpha}})) \to \lambda_{A_{\alpha}}^{A_{\alpha +1}}((I^{\omega})^{A_{\alpha}})\]

\noindent which acts fiberwise (i.e. $\pi_{A_{\alpha}}^{A_{\alpha +1}}g =\pi_{A_{\alpha}}^{A_{\alpha +1}}$). Here is the expression for $g$:

\[ g = \lambda_{A_{\alpha}}^{A_{\alpha +1}}h_{\alpha}p_{A_{\alpha}}^{A_{\alpha +1}}f^{-1}|f(i_{A_{\alpha}}^{A_{\alpha +1}}(Y_{A_{\alpha}})) .\]

By the fibered $Z$-set unknotting theorem \cite{tw}, $g$ extends to a homeomorphism $G \colon (I^{\omega})^{A_{\alpha +1}} \to (I^{\omega})^{A_{\alpha +1}}$ such that $\pi_{A_{\alpha}}^{A_{\alpha +1}} G = \pi_{A_{\alpha}}^{A_{\alpha +1}}$ (i.e. $G$ acts fiberwise). Then the required homeomorphism $h_{A_{\alpha +1}}$ is defined as the composition $Gf \colon Y_{A_{\alpha +1}} \to (I^{\omega})^{A_{\alpha +1}}$. Straightforward verification shows that all the needed properties are satisfied.

This completes the inductive process. Now let $h = \lim\{ h_{\alpha} \colon \alpha < \omega_{1}\}$. It is easy to see that $h \colon Y \to I^{A}$ is a homeomorphim such that $h(i_{A_{\alpha}}(Y_{A_{\alpha}})) = \lambda_{A_{\alpha}}((I^{\omega})^{A_{\alpha}})$ for each $\alpha < \omega_{1}$. Consequently, $h(X) = \Sigma$.
\end{proof}

\begin{cor}\label{C:retract}
Let $X$ be a retract of $\Sigma$. Then $X \times \Sigma$ is homeomorphic to $\Sigma$.
\end{cor}
\begin{proof}
Note that $X \times \Sigma$ is a retract of $\Sigma \times \Sigma \approx \Sigma$ and has no $G_{\delta}$-points.
\end{proof}

\begin{cor}\label{C:productwithcompactum}
The following conditions are equivalent for a compact space $X$:
\begin{itemize}
\item[(i)]
$X \times \Sigma$ is homeomorphic to $\Sigma$,
\item[(ii)]
$X$ is a metrizable absolute retract.
\end{itemize}
\end{cor}
\begin{proof}
(i) $\Rightarrow$ (ii).
Let $h \colon \Sigma \to X \times \Sigma$ be a homeomorphism and $\pi \colon X \times \Sigma \to X$ be the projection. Clearly $r = \pi h \colon \Sigma \to X$ is a retraction. Since $I^{\omega_{1}}$ is the Stone-\v{C}ech compactification of $\Sigma$ (and since $X$ is compact), $r$ admits the extension $\tilde{r} \colon I^{\omega_{1}} \to X$. Therefore $X$, as a retract of $I^{\omega_{1}}$, is an absolute retract. Note also that $X$ is separable (as an image of $I^{\omega_{1}}$). But separable compact subspaces of $\Sigma$ are metrizable.

(ii) $\Rightarrow$ (i). Apply Corollary \ref{C:retract}.
\end{proof}

\end{document}